\documentclass{article}
\begin{document}

Tsemo Aristide

College Boreal

1 Yonge Street, Toronto

tsemo58@yahoo.ca

\bigskip

                               \centerline{\bf Decomposition of groups and top couples.}
\bigskip
\bigskip

\centerline{\bf Abstract.}

{\it Recently, we have endowed various categories of groups with topologies. The purpose of this paper is to introduce
on these categories others topologies which are statistically more suitable to study well-known problems in groups theory.
We use this framework to define a notion of prime ideal and to provide a decomposition of a large class of  groups into a product of prime
Remark that a similar question has been studied in [5] by Kurata with innocent methods. We remark that these topologies can be extended to other categories like the categories of commutative algebras, associative algebras and left symmetric.}

\bigskip

{\bf Definition 1.}

A top couple $(C,D)$ is defined by:

A subcategory $C$ of the category of groups, a subclass $D$ of $C$ which satisfies the
following properties:

T1. Let $G, G'$ be objects of $C$ such that $G'$ is in $D$, if there exists
an injective morphism $i:G\rightarrow G'$, then $G$ is in $D$.

T2. Let $G$ be an object of $D$, $I$ and $J$ two normal subgroups of $G$ such that
 $I\cap J = 1$, then $I=1$ or $J=1$.

T3. Let $G$ be an object of $C$, the normal subgroup $I$ of $G$ is an ideal of $G$ if and
only if the quotient $G/I$ is an object of $C$; we suppose that the inverse image of an ideal by a morphism
of $C$  is an ideal.

\medskip

{\bf Remark.}

\medskip

Let $[I,J]$ be the subgroup generated by the commutators $[x,y]=xyx^{-1}y^{-1}$.
$x\in I$ and $y\in J$.In [12], we have defined a notion of Top couple where we have replaced the axiom T2 by the axiom T'2 as follows:
Let $G$ be an object of $D$, $I$, $J$ two normal subgroups of $G$, $[I,J]=1$ implies $I=1$ or $J=1$;
remark that $[I,J]\subset I\cap J$. This enables  to obtain more examples
of Top couples which are eventually commutative and non trivial. We start by our first example:

\medskip

Let $G$ be a group, we denote by $C_G$ the comma category over $G$, the objects
of $C_G$ are morphisms $f_H:G\rightarrow H$. We denote such an object by $(f_H,H)$.
A morphism between $(f_H,H)$ and $(f_L,L)$ is a morphism of groups $f:H\rightarrow L$
such that $f\circ f_H=f_L$.
Let $(H,\phi_H)$ be an object of $C(G)$ and $x$ an element of $H$, we
denote by $G(x)$ the subgroup of $H$ generated by $\{gxg^{-1},
g\in G\}$.  A non trivial element $x$ of $H$
 is a divisor of zero if and only if there exists a non trivial element $y$ of $H$ such that
$G(x)\cap G(y)=\{1\}$  and $[G(x),G(y)]=1$. We denote by $D_G$ the subcategory of $C_G$ whose objects are the objects
of $C_G$ without divisors of zero.

\medskip

{\bf Proposition 1.}
{\it The couple $(C_G,D_G)$  is a Top couple.}

\medskip

{\bf Proof.}
Let us verify the property T1: Let $H$ and $H'$ be elements of $C_G$, suppose that $H'$
is an object of $D_G$  and there exists an injective morphism $i:H\rightarrow H'$. If $x,y$ are elements
of $H$ such that $G(x)\cap G(y) =1$ and $[G(x),G(y)]=1$, we also have $G(i(x))\cap G(i(y)) =1$ and $[G(i(x)),G(i(y))]=1$
since $i$ is injective. Since $H'$ does not have divisors of zero, we deduce that $i(x)=1$ or $i(y)=1$.
This implies that $x=1$ or $y=1$ since $i$ is injective.

The verification of T2:

Let $H$ be an object $D_G$, $I$ and $J$ two normal subgroups of $H$ such that  $I\cap J=1$.
Suppose that $I$ and $J$ are not trivial.
Let $x$ be a non trivial element of $I$ and $y$ be a non trivial element of $J$, we have $G(x)\subset I$ and $G(y)\subset G$,
this implies that $G(x)\cap G(y)\subset I\cap J=1$ and $[G(x),G(y)]\subset [I,J]\subset I\cap J =1$. Since $H$ does not have divisors
of zero, we deduce that $x=1$ or $y=1$. This is a contradiction.

Verification of T3:

Let $f:H\rightarrow H'$ be a morphism of $C_G$, and $I$ an ideal of $H'$; $f^{-1}(I)$ is an ideal of $H$
since we can endow $H/f^{-1}(I)$ with the structure induced by the morphism $p\circ f_H$, where $p:H\rightarrow H/f^{-1}(I)$
is the canonical projection.

\medskip

{\bf Definitions 2.}
 Let $(C,D)$ be a Top couple, and $H$ an object of $C$, an ideal $P$ of $H$ is prime
 if and only if $H/P$ is an object of $D$.

 For every normal subgroup $I$ of $H$, we denote by $V_H(I)$ the set of prime ideals which contain $I$.

 \medskip

 {\bf Proposition 2.}
 {\it Let $(C,D)$ be a Top couple and $H$ an object of $C$. For every normal subgroups $I,J$ of $H$,
 we have $V_H(I\cap J)=V_H(I)\bigcup V_H(J)$.

 Let $(I_a)_{a\in A}$ be a family of normal subgroups of $H$, and $I_A$ the normal subgroup generated by $(I_a)_{a\in A}$,
 we have $V_H(I_A)=\cap_{a\in A}V(I_a)$.
 }

 \medskip

 {\bf Proof.}
 Firstly, we show that $V_H(I\cap J) = V_H(I)\bigcup V_H(J)$. Let $P$ be an element of $V_H(I\cap J)$ suppose that
 $P$ does not contain neither $I$ nor $J$. Let $x\in I$, $y\in J$ which are not elements of $P$. We denote by $u(x)$
 the normal subgroup of $H$ generated by $x$. We have $u(x)\cap u(y)\subset I\cap J\subset P$.
  This implies that $x\in P$ or $y\in P$.

 Let $P\in V_H(I_A)$. For every $a\in A$, $I_a\subset I_A\subset P$. This implies that $P\subset \cap_{a\in A}V_H(I_a)$.
 Let $P\in \cap_{a\in A} V_H(I_a)$, for every $a\in A$, $I_a\subset P$; this implies that $I_A\subset P$.

\medskip

{\bf Remark.}

\medskip

The space $Spec_G(H)$ of prime ideals is endowed with a topology whose closed subsets are  the subsets $V_H(I)$ and the empty subset of $H$.

Let $x$ and $y$ be divisors of zero in the $G$-group $H$; the subgroup
of $H$ generates by $G(x)$ and $G(y)$ is isomorphic to the direct product $G(x)\times G(y)$.

This leads to the following definitions:

\medskip

{\bf Definitions 3.}
Let $H$ be an element of $C_G$, the adjoint representation $Ad:G\rightarrow Aut(H)$ is the morphism
which associates to $g\in G$ the automorphism of $H$ defined by $Ad(g)(h)=ghg^{-1}, h\in H$.

Let $H$ be an object of $C_G$ a non trivial subgroup $H'$ of $H$ stable by the adjoint representation
 is $G$-decomposable if and only if there exists two non trivial subgroups $H_1$ and $H_2$
stable by the adjoint representations
and an isomorphism of groups $f:H\rightarrow H_1\times H_2$ which commutes with the adjoint representation.

An object $H$ of $C_G$ is locally $G$-indecomposable if every non trivial subgroup of $H$ is not $G$-decomposable.

If $G$ is the trivial group, we will omit the suffix $G$ in the previous definitions, for example,
we will speak of decomposable groups and locally indecomposable groups.

\medskip

{\bf Proposition 3.}
{\it A $G$-group $H$ does not have divisors of zero if and only if $H$ is locally $G$-indecomposable.}

\medskip

{\bf Proof.}
Suppose that the $G$-group $H$ does not have divisors of zero, let $L$ be a subgroup stable by the adjoint action; suppose that
 $L$ is isomorphic to the product of the non trivial subgroups $L_1$ and $L_2$ stable
by the adjoint representation. Let $x_1\in L_1$ and
$x_2\in L_2$ be non trivial elements; $(x_1,1)$ and $(1,x_2)$ are divisors of zero. This is a contradiction.

Conversely, suppose that the $G$-group $H$ is locally indecomposable; let $x$ and $y$ be divisors of zero; the
subgroup of $H$ generates by $G(x)$ and $G(y)$ is a subgroup of $G$ which is the direct product of the subgroups
$G(x)$ and $G(y)$ which are stable by the adjoint action. This is a contradiction.

\medskip

{\bf Remark.}

Let $G$ be a group, to study the geometry of objects of $C_G$, it is very important to know objects
without divisors of zero. Firstly, we are going to study these objects for $G=1$. We are also going to classify
finitely generated nilpotent groups who do not have divisors of zero. Remark that finite groups without
divisors of zero have been classified by Marin when $G=1$; to present his result, let us recall that
the quaternionic group $Q_n$ ($n$ is an integer superior or equal to $3$) is a finite group of order $2^n$
with the presentation:

$$
<x,y :x^{2^{n-1}}=1, y^2 = x^{2^{n-2}}, y^{-1}xy=x^{-1}>
$$

\medskip

{\bf Theorem Marin [6].}
{\it Suppose that $G=1$; a finite group $H$ is indecomposable if and only if:

1. $H$ is isomorphic to $Z/p^n$ for some prime $p$.

2. $H$ is isomorphic to $Q_n, n\geq 3$.

3.  $H$ is isomorphic to an extension of $Z/q^b$ by $Z/p^a$ where $p$ and $q$ are different prime
integers such that $p$ is odd, $q^b$ divides $p-1$ and the image of $Z/q^b$ in $(Z/p^a)^*$ has order $q^b$.

}

\medskip.

{\bf Proposition 4.}
{\it Suppose that $G=1$, let $H$ be a group without divisors of
zero. The rank of every commutative subgroup of $H$ is inferior to $1$. In particular
the rank of the center $C(H)$ is inferior to $1$.  If the center is not trivial, for every $y\in H$, there exists $n\in N$ such
that $y^n$ is an element of $C(H)$ distinct of the identity.
If the order of the center $C(H)$ is finite, then the order of every element of
$H$ is finite and in this case the order of such an element is $p^n$ where $p$ is a prime integer.}

\medskip

{\bf Proof.}
If the rank of a commutative subgroup $L$ is strictly greater than $1$, there exists
non trivial elements $x,y$ in $L$ such that $[x,y]=1$ and
$(x)\cap (y)=1$. Where $(x)$ is the subgroup of $H$ generated by $x$. This is in contradiction with the fact that $H$
does not have zero divisors. Let $z\in C(H)$ be a non trivial element, for every
element $x \in H$, we have $[x,z]=1$, since $H$ does not have
divisors of zero, we deduce that $(x)\cap (z)$ is not the trivial
group.

Suppose that the center of $H$ has a finite order, for any element $x\in H$,
there exists an integer $n$ such that $x^n\in C(H)$, $x^n$ and henceforth $x$ has a finite order.
If the order of $z$ is the product $nm$ of two integers $n$ and $m$ which are relatively prime, then
$z^n$ and $z^m$ are divisors of zero.

\medskip

{\bf Theorem 1.}
{\it Suppose that $G=1$, let $H$ be a finitely generated nilpotent group
without divisors of zero. Then $H$ is finite or $H$ is isomorphic to $Z$.}

\medskip

{\bf Proof.}
Let $H$ be a non trivial finitely generated nilpotent group.
 recall that the derivative sequence of $H$
is defined by $H^0=H$, and $H^{(n)}=[H,H^{(n-1)}]$. There exists $n$ such that $H^{(n)}=1$,
and $H^{(n-1)}$ is not trivial and contained in the center of $H$. The proposition 4 shows that the
rank of $H^{(n-1)}$ is $1$. Suppose that there exists an element $x$ of $H^{(n-1)}$ which
has a finite order, then every element of $H$ has a finite order. The subgroup $H^{(n-2)}$ is
finite since it is the extension of a commutative finite group by a commutative finite group;
recursively, we obtain that $H$ is finite.

Suppose  now that $C(H)$ has infinite order and the rank of $H$ is different of $1$.
We have $[H,H^{(n-2)}]=H^{(n-1)}$. This implies the existence of an element $x\in H$ and $y\in H^{(n-2)}$
such that $[x,y]\in H^{(n-1)}$ and is distinct of the neutral element and has an infinite order. Remark that $[x,y]=h$
is in the center of $H$. There exists integers $n,m$ such that $x^n\in C(H)$ and $y^m\in C(H)$.
We have $x^ny^mx^{-n}y^{-m}=h^{mn}y^mx^nx^{-n}y^{-m}=1$. This implies that the order of $h$ is finite.
This is a contradiction with the hypothesis.

\medskip

{\bf Corollary 1.}
{\it A finitely generated locally indecomposable whose commutator subgroup is nilpotent is a finite group or is a finite extension of $Z$.}

\medskip

{\bf Proof.}
Let $H$ be a finitely generated locally indecomposable whose commutator subgroup is nilpotent. Then $[H,H]$ is a locally indecomposable
nilpotent group. Suppose that $[H,H]$ is infinite, thus $[H,H]=Z$. Let $x$ be an element of $H$; $Ad_x:[H,H]\rightarrow [H,H]$ defined by $Ad_x(y)=xyx^{-1}$
has order inferior to $2$ since the group of automorphisms of $Z$ is isomorphic to $Z/2$. Let $y$ be a
non trivial element of $[H,H]$, we deduce that for every $x\in H$, $[x^2,y]=1$. Since $H$ does not have
divisors of zero, it results that there exists $n,m$ such that $x^{2n}=y^m$. Thus the quotient $H/[H,H]$
is finite.

Suppose that $[H,H]$ is finite and for every $x\in H$, $Ad_x$ is an automorphism of a finite group, thus there
exists $n$ such that $Ad_{x^n}$ is the identity. Let  $z$ be a non trivial element of $[H,H]$,
 $[x^n,z]=1$, since $H$ does not have divisors of zero, we deduce that there exists $m$ such that $x^{nm}\in [H,H]$;
 thus every element of $H$ has a finite order. Since $H$ is solvable, we deduce that $H$ is finite.

\medskip

{\bf Corollary 2.}
{\it A subgroup $I$ of a finitely generated commutative group $H$ is a prime ideal if and only if
$G/H$ is isomorphic either to $Z$ or to $Z/p^n$ where $p$ is a prime.}

\medskip

{\bf Proof.}
Let $I$ be a prime ideal of the finitely generated commutative group $H$, if $H/I$ is finite, Marin implies that
$H/I$ is isomorphic to $Z/p^n$ where $n$ is a prime if $H/I$ is infinite, since it nilpotent, proposition implies that
$H$ is isomorphic to $Z$.

\medskip

{\bf Remark.}

Suppose that $H=Z$ the group of relative integers. Let $I$ be a
ideal of $H$, we know that $I$ is a subgroup generated by a
positive integer $n$, write $n=\prod_{i\in I}p_i^{n_i}$. Let $p$
be a prime number  and $a$ and integer, the prime ideal $(p^a)$
generated by $p^a$ is an element of $V((n))$ if and only if $p^a$
divides $n$.

We are going to present other examples of locally indecomposable groups
Recall that the Tarski group is an infinite group $H$ such that there exists a prime integer $p$ such that every
subgroup of $H$ is isomorphic to the cyclic group $Z/p$. The Tarski group is known to be simple.
Olshans'kii [8] and have shown the existence of Tarski groups for $p>10^{75}$.

Adyan and Lysenok [1] and have generalized the construction of Ovshan'skii and shown that for $n> 1003$ there
exists non commutative groups $H$ such that every proper subgroup of $H$ is isomorphic to  a subgroup
isomorphic to $Z/n$, we will call these groups Adyan-Lysenok groups.

Remark that the Adyan-Lysenok groups $H$ defined for $n=p^m$  is a domain for $G=1$:
Let  $x$, $y$ divisors of zero in $H$, since the subgroup $<x,y>$ generated by $x$ and $y$ is
a commutative subgroup  we deduce that $<x,y>$
is isomorphic to a subgroup of $Z/p^m$. This is in contradiction with the fact that $<x>\cap <y>$ is trivial.

More domains can be constructed by using the following proposition:

\medskip

{\bf Proposition 5.}
{\it The free product two locally indecomposable groups is a locally indecomposable group.}

\medskip

{\bf Proof.}
Let $G$ and $H$ be two locally indecomposable groups. Let $x$ and $y$ be divisors of zero,
then since $xy=yx$, the corrollary  [7] 4.1.6 p.187 shows either:

- $x$ and $y$ are conjugated in the same factor of $G$ or $H$. This is impossible since
$G$ and $H$ are locally indecomposable

- $x$ and $y$ are the power of the same element. This is in contradiction with the fact that
$x$ and $y$ are divisors of zero.

\medskip

{\bf Definition 4.}
Let $H$ be an element of $C(G)$, we denote by $Rad_G(H)$ the
intersection of all the prime ideals of $H$.

\medskip

Recall that a topological space $X$ is irreducible if and only if
it is not the union of two proper subsets.

We say that an ideal $I$ is a radical ideal if it is the
intersection of all the prime which contains $I$.

\medskip

{\bf Proposition 6.}
{\it Let $H$ be an element of $C(G)$, and $I$ a radical ideal of
$H$, then $V_H(I)$ is irreducible if and only if $I$ is a prime.}

\medskip

{\bf Proof.}
Suppose that $I$ is a prime, and $V_H(I)=V_H(J)\bigcup V_H(K)$ where $V_H(J)$ and
$V_H(K)$ are proper subsets,
since $I$ is a prime, $I$ is an element of $V_H(I)$. This implies
that $I\in V_H(J)$ or $V_H(K)$. If $I$ is an element of $V_H(J)$,
then $V_H(I)\subset V_H(J)$; if $I\in V_H(K)$, then $V_H(I)\subset
V_H(K)$.This is a contradiction with the fact that $V_H(J)$ and
$V_H(K)$ are proper subsets of $V_H(I)$.

Suppose that $V_H(I)$ is irreducible; let $x,y$ be elements of $H$
such that $[G(x),G(y)]\subset I$ and $G(x)\cap G(y)\subset I$. Let
$u(x)$ be the normal subgroup generated by $x$, $u(x)\cap u(y)$
and $[u(x),u(y)]$ are contained in $I$. Since $V_H(u(x)\cap u(y)) = V_H(u(x))\bigcup V_H(u(y))$,
this implies that $V_H(I)= V_H(u(x))\cap V_H(I)\bigcup V_H(u(y))\cap V_H(I)$. Since $V_H(I)$
is irreducible, we deduce that $V_H(I)$
is contained in $V_H(u(x))$ or is contained in $V_H(u(y))$. If
$V_H(I)$ is contained in $V_H(u(x))$, the $\cap_{P\in
V_H(I)}P=I$ contains $u(x)$. It results that $x\in I$ since $I$ is a radical ideal. Similarly,
if $V_H(I)\subset V_H(u(y))$ we deduce that $y\in I$.

\medskip

{\bf Definition 5.}
Recall that a space is Noetherian if and only if every ascending
chain of closed subsets $Z_0\subset Z_1...\subset Z_n\subset...$ stabilizes,
this is equivalent to saying that there exists  $i$ such that for every $n>i, Z_n=Z_i$. We deduce
that the topological space $Spec_G(H)$ is Noetherian if and only if
a descending chain of normal subgroups of $H$ $(I_n)_{n\in N}$ such that
$I_{n+1}\subset I_n$ stabilizes.

\medskip

{\bf Remark.}

Let  $G$ be a group:

 $G$ is an element of $Spec_G(G)$, we denote by $Spec_G(G)^*$, $Spec_G(G)-\{G\}$. We will denote by $V^*_H(I)$,
 the intersection $V_H(I)\cap Spec_G(G)^*$.

For every element $x\in G$, $G(x)$ is the normal subgroup generated by $x$.

A maximal normal subgroup $I$ of $G$ is a prime, since $G/I$ is a simple group.

\medskip

{\bf Theorem 2.}
{\it Suppose that $Spec_G(G)^*$ is Noetherien and $Rad_G(G)=1$, then
$G$ is the product of groups $G_1\times...\times G_n$ such that
 for every $i$, the subgroup $H_i$ of $H$
generated by $G_j, j\in\{1,...,n\}-\{i\}$ is a prime. Moreover, this
decomposition is unique up to the permutation of the $G_i$.}

\medskip

{\bf Proof.}
Suppose that $Spec_G(G)^*$ is Noetherian, then $Spec_G(G)^*$ is the
disjoint union of closed subsets $(V_G(H_i))_{i=1,...,n}$.

The intersection $\cap_{i=1,...,n}H_i=1$ . This is due to the fact
that $V^*_G(\cap_{i=1,...,n}H_i)=V^*_G(H_1)\bigcup...\bigcup V^*_G(H_n)
= Spec^*_G(G)$ and $Rad_G(G)=1$.

We write $G_i=\cap_{j\in\{1,...,n\}-\{i\}} H_j$. We are going to
show that $G$ is isomorphic to the direct product
$G_1\times...\times G_n$.

Firstly, remark that $G_i\cap G_j =\cap_{k=1,..,n}H_k=1$ if $i\neq
j$. Since the subgroup $G_i$ are normal, for $i\neq j$, we have $[G_i,G_j]\subset G_i\cap G_j=1$.
This implies that the subgroup $L$ of $H$ generated by $(G_i)_\{i=1,...,n\}$ is isomorphic to the
direct product $G_1\times G_2\times...\times G_n$.
It remains to shows that $G$ is equal to its subgroup $L$.

We have $V^*_G(G_i)=\bigcup_{j\in \{1,..,n\},j\neq i}V_G^*(H_i)$. This implies
that $V_G^*(L)=\cap_{i=1,...,n}\bigcup_{j\in \{1,..,n\},j\neq i}V^*_G(H_i)$ is empty.
We deduce that $L=H$, otherwise $L$ would have been contained in a maximal ideal which would have been
an element of $V_G^*(L)$.

We show now that the subgroup $L_i$  of generated by $(G_j)_{j\neq i}$ is $H_i$.
 For every $j\neq i$, $G_j\subset H_i$. Suppose that there exists an element $x\in H_i$ which is not in $L_i$.
 Since $H=G_1\times...\times G_n$, we can write $=(x_1,...,x_n), x_j\in G_j$ and $x_i\neq 1$, we have $x_j\in H_i, j\neq i$. This implies
 that $x_i\in H_i$. This is a contradiction since $H_i\cap G_i=\{1\}$

We show now that the decomposition is unique. Suppose that there are two decompositions
$H=G_1\times...\times G_n$ and $H=U_1\times...\times U_m$ such that the group $H_i$ generated by $1\times..G_j\times 1.., j\neq i$
is a prime ideal, the group $L_i$ generated by $1\times..U_j\times..\times 1 j\neq i$ is also a prime ideal.
Then $\bigcup_{i=1,..,n}V^*_G(H_i)$ and $\bigcup_{i=1,..,m}V_G^*(L_i)$ are decomposition
of $Spec_G^*(G)$ as union of irreducible components. Since this decomposition is unique, we deduce that
$n=m$, and up to permutation that $V_G^*(H_i)=V_G^*(L_i)$, since $U_i$ and $H_i$ are prime, we deduce that
$H_i=L_i$. This implies that $G_i\simeq G/H_i$ and $U_i\simeq G/L_i$ are isomorphic.

\medskip

{\bf Corollary 3.}
{\it Suppose that $G$ is a finite group and $Rad_G(G)=1$, then $G$
is a product of indecomposable subgroups.}

\bigskip

{\bf Some generalizations.}

Let $A$ be a commutative ring, in classical algebraic geometry a prime ideal
$P$ of $A$ is an ideal $P$ such for every elements $a,b\in A$, $ab\in P$ implies that
$a\in P$ or $b\in P$. Inspired by the topologies defined above, we define the following notion:

\medskip

{\bf Definitions.}
Let $A$ be a ring non necessarily commutative, $a,b$ elements of $A$. We denote by $I(a)$ the two-sided ideal generated by $a$.
 A two-sided ideal of the ring $A$ is a $p$-prime if for every elements $a,b\in A$, $I(a)\cap I(b)\in P$
implies that $a\in P$ or $b\in P$.

Let $I$ be a two-sided ideal of $A$, we denote by $V(I)$ the set of prime ideals of $A$ which contain $I$.

\bigskip

{\bf Proposition.}
{\it Let $I,J$ be two-sided ideals of $A$, we have: $V(I\cap J)=V(I)\bigcup V(J)$. Let $(I_a)_{a\in A}$ be family
of ideals of $A$ which generates the ideal $I_A$, we have $V(I_A)=\cap_{a\in A}V(I_a)$.}

\bigskip

{\bf Proof.}
Firstly, we show that $V(I\cap J)=V(I)\bigcup V(J)$. Since $I\cap J\subset I$ and $I\cap J\subset J$,
we have $V(I)\subset V(I\cap J)$ and $V(J)\subset V(I\cap J)$. Let $P$ be an element of $V(I\cap J)$,
suppose that $P$ does not contain $I$ and $J$. Let $a\in I, b\in J$ be elements which are not in $P$, $I(a)\cap I(b)\subset I\cap J$.
 This is a contradiction since $P$ is a prime ideal.

Let $P$ be an element of $V(I_A)$, since $I_A\subset P$, $I_a\subset P$ for every $a\in A$,
this implies that $P\in \cap_{a\in A}V(I_a)$. Conversely, let $P\in \cap_{a\in A}V(I_a)$,
for every $a\in A$, $I_a\subset P$. This implies that $I_A\subset P$.

\medskip

{\bf Examples.}

\bigskip

Suppose that $A$ is a commutative algebra, an ideal $I$ is a prime if and only if
for every $a,b\in A$, $I(a)\cap I(b)\subset P$ implies that $a\in P$ or $b\in P$.
This structure is different from the classical notion of prime. As we have seen, if $A=Z$, $Z/p^n$
is a prime ideal.

\bigskip
\bigskip
\bigskip

{\bf References.}

\smallskip

 1. Adyan S.I. Lysenok I.G., "On groups all of whose proper subgroups of which are finite cyclic", Izv. AN SSSR. Ser. matem., 55:5 (1991), 933–990

2. Amaglobeli. M.G Algebraic sets and coordinate groups for a free nilpotent group of nilpotency class 2. Sibirsk. Mat. Zh. Volume 48 p. 5-10.

 3.  Baumslag, G,  Miasnikov, A.  Remeslennikov, V.N. Algebraic geometry over groups I. Algebraic sets and ideal theory. J. Algebra. 1999, 219, 16–79.

 4.. A.Grothendieck, \'El'ements de g\'eom\'etrie alg\'ebrique I.Publications math\'ematiques de l'I.H.E.S 4, 5-228

5. Kurata, Y. A decomposition of normal subgroup in a group. Osaka J. Math.
1 (1964), 201-229

6. Marin, I. Strongly indecomposable finite groups Expositiones Mathematicae  , vol. 26, no. 3, pp. 261-267, 2008

7. Magnus, Karass, Solitar.  Combinatorial group theory, Dover publication 1976

 8. Olshanskii A., Groups of bounded period with subgroups of prime order, Algebra and Logic 21 (1983), 369-418; translation of Algebra i Logika 21 (1982), 553-618

9. Scott W.R Algebraically closed groups Proc. Amer. Math. Soc. 2 (1951) 118-121

10. Serre J-P. (GAGA) G\'eom\'etrie alg\'ebrique, g\'eom\'etrie analytique, Annales de l'Institut Fourier, Grenoble
 t. 6 1955-1956 1-42.

 11. Tsemo, A. Scheme theory for groups and Lie algebra, International Journal of Algebra 5. 2011 139-148

12. Tsemo, A. Some properties of schemes in groups theory and Top couples International Journal of Algebra, Vol. 7, 2013, no. 1, 25 - 48.

13. Tsemo, A. Theory of curves, in preparation.

\end{document}